\documentclass[12pt]{amsart}
\usepackage{amsmath,amsthm,latexsym,amscd,amsbsy,amssymb,url}
\usepackage[all]{xypic}
 \relax

\newtheorem{thm}{Theorem}[section]

\newtheorem{ex}[thm]{Example}
\newtheorem{que}[thm]{Question}
\newtheorem{problem}[thm]{Problem}

\begin{document}


\title
{Better Subtopologies}

\author{Alexander Arhangel'skii}
\email{arhangel.alex@gmail.com}
\address{Moscow, Russia}

\author{Raushan Buzyakova}
\email{Raushan\_Buzyakova@yahoo.com}
\address{Department of Mathematics, Polk State College, USA}

\keywords{continuous bijection, continuous injection, linearly ordered space, generalized ordered space, space of ordinals}
\subjclass{54E99,54F05}


\begin{abstract}{
We study conditions under which a space that has a good property and a courser topology with another good property admits a continuous bijection onto a space with both properties.
}
\end{abstract}

\maketitle
\markboth{A. Arhangelskii and R. Buzyakova}{Better Subtopologies}
{ }

\section{Introduction}\label{S:introduction}

\par\bigskip\noindent
This note is an invitation to a discussion around the following general problem:
\par\bigskip\noindent
\begin{problem}\label{problem:p1}
Let $X$ have a property $A$ and admit a continuous bijection onto a space with property $B$. Under what additional conditions, does $X$ admit a continuous bijection onto a space with both properties $A$ and $B$?
\end{problem}

\par\bigskip\noindent
The next problem is a natural generalization of Problem \ref{problem:p1}.

\par\bigskip\noindent
\begin{problem}\label{problem:p2}
Let $X$ admit a continuous bijection onto a space with property $A$ and onto a space with property $B$. Under what additional conditions, does $X$ admit a continuous bijection onto a space with both properties $A$ and $B$?
\end{problem}

\par\bigskip\noindent
In search for affirmative results it makes sense to consider very restrictive properties $A$ and $B$ such as 
\begin{enumerate}
    \item $A$ is a metrizable space and $B$ is a topological group.
    \item $A$ is metrizable and $B$ is orderable.
    \item $A$ is separable and metrizable and $B$ is compact.
    \item $A$ is metrizable and $B$ is homogeneous.
    \item any of the above combinations in the reverse order.
\end{enumerate}
A notable result that gives hope for affirmative outcomes in the direction of this list belongs to Pasynkov \cite{P}, which states, in particular,  that any separable submetrizable space of Ind-dimension at most $n$ admits a continuous bijection onto a metrizable space of Ind-dimension at most $n$. This Pasynkov's theorem is a generalization of the classical factorization theorem of Mardesic \cite{M}.  

\par\bigskip\noindent
In this study, we will present some affirmative results and support their relative generality  by relevant examples. We will identify some more specific questions related to the above problems that may pave routes for future research.

\par\bigskip\noindent
Note that for a topological space,   admitting a continuous bijection onto a space with a desired property is equivalent to having a coarser topology (that is, a subtopology) with that property. We will, therefore, use these three terms interchangeably. 

\par\bigskip\noindent
In notation and terminology, we will follow \cite{Eng}. In our study, we will concentrate on Problem \ref{problem:p1} with $A$ or $B$ or both being (nice) ordered spaces or their subspaces. Recall that a linearly ordered topological space (abbreviated as LOTS) is an ordered set $\langle L, <\rangle$ with at  least two points and a topology generated by sets in form $\{x\in L: x<a\}$, $\{x\in L: x>a\}$. Of course, a singleton is also a LOTS but is of no interest to us. A generalized ordered space, abbreviated as a GO-space,  is any subspace of a LOTS. It is the result of \v Cech (for a proof, see, for example \cite{BL}) that a Hausdorff topological space $X$ with a linear order $<$ is a generalized ordered space if and only if the topology of $X$ is generated by a collection of $<$-convex subsets of $X$. To distinguish an open interval from an ordered pair, the latter will be denoted by $\langle a,b\rangle$. The lexicographical product of linearly ordered sets $L$ and $M$ is denoted by $L\times_l M$.

\par\bigskip\noindent
\section{Study}

\par\bigskip\noindent
We begin our study by an affirmative result that justifies the rest of our discussion.
\par\bigskip\noindent
\begin{thm}\label{thm:go1}
Let $X$ be a first-countable space that admits a continuous bijection onto a GO-space $Y$. Then $X$ admits a continuous bijection onto a first-countable GO-space.
\end{thm} 
\begin{proof} Fix a continuous bijection $f$ of $X$ onto $Y$. We will construct a first-countable GO-space $Z$ by refining the topology of $Y$. For this we need the following statement.
\par\bigskip\noindent
{\it Claim. Let $y\in Y$.  Suppose that  $x=f^{-1}(y)$ is  a limit point for $A=f^{-1}(\{z:z<y\})$. Then $y$ has countable  character in $\{z:z\leq y\}$.
}
\par\smallskip\noindent
To prove the claim, we note that,  by first countability there exists a sequence $\{x_n\}$ of elements of $A$ that converges to $x$. Then, the sequence $\{f(x_n)\}$ consists of elements of $\{z: z<y\}$ and converges to $y=f(x)$, which completes the proof of the claim. 
\par\bigskip\noindent
\underline{Construction of $Z$}: The underlying set of $Z$ is that of $Y$. A basis ${\mathcal B}_Z$ for the topology of $Z$ consists of sets that fall in one of the following three categories:
\begin{enumerate}
	\item If $I$ is a convex open set in $Y$, then $I\in {\mathcal B}_Z$.
	\item If $y$ has uncountable character in  $\{z:z\leq y\}$, then  $\{z:z\geq y\}\in {\mathcal B}_Z$.
	\item If $y$ has uncountable character in  $\{z:z\geq y\}$, then  $\{z:z\leq y\}\in {\mathcal B}_Z$.
\end{enumerate}
Due to (1), the identity map  from $Z$ to $Y$ is a continuous bijection. Since ${\mathcal B}_Z$ consists of convex sets of $Y$, the space $Z$ is a GO-space (see the introduction section for reference). Due to (2) and (3), $Z$ is first-countable. Due to {\it Claim}, the map $g:X\to Z$ defined by $g(x)=f(x)$ is a continuous bijection with respect to the topology of $Z$ as well.

\par\bigskip\noindent
Our proof is complete. For further reference we would like to make the following final remark.
\par\bigskip\noindent
{\it Remark:
Note that if $Y$ is a subspace of an ordinal, then $Z$ is homeomorphic to a subspace of an ordinal too. 
}
\end{proof}

\par\bigskip\noindent
It is natural
to wonder if one can replace each occurrence of "GO-space" by "LOTS" in the statement of Theorem \ref{thm:go1}. The following example shows that the answer is negative if no additional assumptions are made.

\par\bigskip\noindent
\begin{ex}\label{ex:go22} 
There exists a first-countable GO-space $L$ that admits a continuous bijection onto a LOTS but not onto a first-countable LOTS.
\end{ex}
\begin{proof} Let $I$ be the single-point compactification of  $\omega_1\times_l [0,1)$ by adjoining the point $\infty$. In literature, this space is often called "the long segment" (see, for example \cite{Eng}). 
The underlying set for our example is  $\mathbb R\times_l I$. Let $\mathcal T_O$ be the order topology on this set. Our space $L$ is the
set $\mathbb R\times_l I$ endowed with the topology generated by 
$\mathcal T_O\cup \{\{\langle x, \infty \rangle \}: x\in \mathbb R\}$, that is, we take the original linearly ordered space 
$\mathbb R\times_l I$ and declare all points of uncountable character isolated. This action makes the newly formed topology first-countable. Since the topology of $L$ has a basis consisting of convex sets, $L$ is a GO-space. The space $L$ admits a continuous bijection onto linearly ordered space  $\mathbb R\times_l I$. Next, let us show that $L$ does not admit a continuous bijection onto a first-countable LOTS.
\par\medskip\noindent
Assume the contrary and let $f$ be a continuous bijection of $L$ onto a first-countable LOTS $M$. For each $x\in \mathbb R$, put $T_x=\{x\}\times_l (I\setminus \{\infty\})$. Due to connectedness of $T_x$, the following statement holds.
\par\bigskip\noindent
{\it Claim 1. $f|_{T_x}$ is order-preserving or order-reversing for each $x\in \mathbb R$.
}
\par\bigskip\noindent
Without loss of generality, we may assume that $f|_{T_x}$ is order-preserving for uncountably many $x\in \mathbb R$. Denote the set of all such $x$ by $A$. Fix $x\in A$. Since $M$ is a first-countable LOTS, $\sup f(T_x)$ does not exist. 
\par\bigskip\noindent
{\it Claim 2. Let $z< f(x,0)$, where $x\in A$. Then there exists a rational number $q$ such that $z<f(q,0)<f(x,0)$.
}
\par\smallskip\noindent
To prove the claim select a strictly increasing sequence $\{q_n\}_n$ of rational numbers that converges to $x$. Then, $\{f(q_n,0)\}_n$ converges to $f(x,0)$.

\par\bigskip\noindent
By Claim 2, we can find distinct $p,q\in \mathbb Q$ such that $f( x, \infty)$ and $f(x,0)$ are separated by the interval $[f(q,0), f(p,0)]$. Without loss of generality, there exists an uncountable $B\subset A$ such that $f(x,\infty) < f(q,0)$ and $f(x,0)>f(p,0)$ for all $x\in B$. Let $\{b_n\}$ be a strictly increasing sequence of elements of $B$ converging to some $b$. Then $\{\langle b_n, 0\rangle\}$ and $\{\langle b_n,\infty\rangle\}$ converge to $\langle b, 0\rangle$. However, by the choice of $B$, the limit of $\{f(b_n,\infty )\}$ and the limit of $\{f(b_n,0)\}$ are separated by the interval $[f(q,0), f(p,0)]$. This contradicts continuity of $f$.
\end{proof}

\par\bigskip\noindent
In view of Theorem \ref{thm:go1} and Example \ref{ex:go22}, it is worth noting that We can replace each occurrence of "GO-space" with "subspace of an ordinal" in  the statement of Theorem \ref{thm:go1} without jeopardizing the validity of the statement.

\par\bigskip\noindent
\begin{thm}\label{thm:go3}
Let $X$ be a first-countable space that admits a continuous bijection onto a subspace of an ordinal. Then $X$ admits a continuous bijection onto a first-countable subspace of an ordinal. 
\end{thm}
\begin{proof}
The proof is similar to that of Theorem \ref{thm:go1} and the needed additional step is elaborated in the end-of-proof remark of Theorem \ref{thm:go1}.
\end{proof}

\par\bigskip\noindent
It is natural to wonder if the first-countability requirement in Theorem \ref{thm:go3} is important. This justifies  the following question.
\par\bigskip\noindent
\begin{que}\label{que:go3}
Does there exist a space of countable pseudocharacter that admits a continuous bijection onto a subspace of an ordinal  but does not admit a continuous bijection onto a first-countable subspace of an ordinal?
\end{que}

\par\bigskip\noindent
Replacing countable character with countable pseudocharacter in Theorem \ref{thm:go3} does lead to a positive outcome under additional conditions.
\par\bigskip\noindent
\begin{thm}\label{thm:go5}
Let $X$ have countable pseudocharacter and admit a continuous bijection onto a hereditarily paracompact subspace of an ordinal. If $ind(X) = 0$, then $X$ admits a continuous bijection onto a first-countable hereditarily paracompact  subspace of an ordinal. 
\end{thm}
\begin{proof}
We will argue by induction. Suppose that for every ordinal $\beta<\alpha$ it is true that if a space of countable pseudocharacter and ind-dimension 0 continuously injects onto a hereditarily paracompact subspace of  $\beta$ then the space continuously bijects onto a first-countable hereditarily paracompact subspace of an ordinal. 
We may assume that $\alpha$ is infinite. We will consider the only three possible cases.

\begin{description}
        \item[\rm Case($\alpha$  is limit)] Since $f(X)$ is paracompact, it can be written as the union of a discrete collection of convex clopen sets each of which is a subset of some ordinal strictly less than $\alpha$. Applying the assumption to each set in the collection, the conclusion follows.
	\item[\rm Case ($\alpha = \beta+1$ and $\beta$ is limit)] Due to the previous case, we may assume that $\beta$ is in the range of $f$. Let $z = f^{-1}(\beta)$. By countable pseudocharacter and zero ind-dimentionality, there exists a disjoint family $\{A_n\}$ of clopen subsets  of $X$ such that $X\setminus \{z\} = \cup \{A_n : n \in \omega\}$. For each $A_n$, $f(A_n)$ is a heridatarily paracompact subspace of $\beta$. By assumption, we can fix a continuous bijection $f_n$ of $A_n$ onto a first-countable hereditarily paracompact subspace of some ordinal. Let $L$ be the GO-space defined as follows: $L = (\oplus_n f(A_n)) \cup \{z\}$. A basis at $z$ consists of intersections of finitely many sets in form $L\setminus f(A_n)$.  Clearly, $L$ is homeomorphic to a first-countable hereditarily paracompact subspace of an ordinal and is a continuous image of $X$.
	\item[\rm Case($\alpha = \beta+1$ and $\beta$ is isolated)] Then $\alpha$ and $\beta$ are topologically homeomorphic, and therefore, the conclusion follows. 
\end{description}
\end{proof}

\par\bigskip\noindent
Clearly, there are too many assumptions in Theorem \ref{thm:go5}. Can we drop some?

\par\bigskip\noindent
\begin{que}
Let $X$ have countable pseudocharacter. Suppose $X$ admits a continuous bijection onto a subspace of an ordinal. Does $X$ admit a continuous bijection onto a first-countable subspace of an ordinal? What if $X$ has  dimension (ind, Ind, or dim) zero?
\end{que}

\par\bigskip\noindent
We would like to finish our study with a few specific questions stemming from Problems \ref{problem:p1} and \ref{problem:p2} that detour from the class of ordered spaces and their subspaces.
\begin{que}
Let $X$ be a locally compact space that admits a continuous bijection onto a separable space. Does $X$ admit a continuous bijection onto a separable compact space.
\end{que}

\par\bigskip\noindent
\begin{que}
Let $X$ be a topological group that admits a continuous bijection onto a compact (paracompact, Chech-complete) space. Does X admit a continuous bijection onto a compact (paracompact, Chech-complete) topological group? Is $X$ is paracompact topological group?
\end{que}

\par

\end{document}